\documentclass[12pt]{article}
\usepackage{amsfonts, amsmath, latexsym, epsfig}
\usepackage{epsf}
\usepackage{url}

\newtheorem{lem}{Lemma}
\newtheorem{cor}{Corollary}
\newcommand{\qed}{\hfill $\Box$ }
\newcommand{\proof}{\noindent{\bf Proof.}\ \ }

\begin{document}
\newcommand{\RR}{\ensuremath{\mathbb{R}}}
\newcommand{\NN}{\ensuremath{\mathbb{N}}}
\newcommand{\QQ}{\ensuremath{\mathbb{Q}}}
\newcommand{\CC}{\ensuremath{\mathbb{C}}}
\newcommand{\ZZ}{\ensuremath{\mathbb{Z}}}
\newcommand{\KK}{\ensuremath{\mathbb{K}}}
\newcommand{\TT}{\ensuremath{\mathbb{T}}}

\author{Mathieu DUTOUR SIKIRI\'C\\
        \normalsize  Institut Rudjer Boskovi\'c, Zagreb\\
\and
        Viatcheslav GRISHUKHIN\\
        \normalsize CEMI RAN, Russia
}

\title{How to compute the rank of a Delaunay polytope}
\date{}

\maketitle

\begin{abstract}
Roughly speaking, the rank of a Delaunay polytope (first introduced in 
\cite{DGL92}) is its number of degrees of freedom.
In \cite{DL}, a method for computing the rank of a Delaunay polytope $P$
using the hypermetrics related to $P$ is given.
Here a simpler more efficient method, which uses affine dependencies instead
of hypermetrics is given.
This method is applied to classical Delaunay polytopes.

Then, we give an example of a Delaunay polytope, which does not have
any affine basis.
\end{abstract}

\section{Introduction}
A lattice $L$ is a set of the form $v_1\ZZ+\dots +v_n\ZZ\subset \RR^n$.
A {\em Delaunay polytope} $P$ is inscribed into an {\em empty sphere}
$S$ such that no point of $L$ is inside $S$ and the vertex-set of
$P$ is $L\cap S$. The Delaunay polytopes of $L$ form a partition
of $\RR^n$.

The vertex-set $V=V(P)$ of a Delaunay polytope $P$ is the support
of a distance space $(V,d_P)$. For $u,v\in V(P)$, the distance
$d_P(u,v)=\Vert u-v\Vert^2$ is the {\em Euclidean norm}
of the vector $u-v$.
A {\em distance vector} $(d(v,v')  )$ with $v,v'\in V$ is called a 
{\em hypermetric} on the set $V$ if it satisfies
$d(v,v')=d(v',v)$, $d(v,v)=0$ and the following {\em hypermetric inequalities}:
\begin{equation}\label{Definition-of-hypermetrics}
H(b)d=\sum_{v,v'\in V} b_vb_{v'}d(v,v')\leq 0\mbox{~for~any~}
b=(b_{v})\in \ZZ^{V}
\mbox{~with~}\sum_{v\in V}b_v=1\,\,.
\end{equation}
The set of distance vectors, satisfying (\ref{Definition-of-hypermetrics}) is called the {\em hypermetric cone} and denoted by $HYP(V)$.

The distance $d_P$ is a hypermetric, i.e., it belongs to the hypermetric
cone $HYP(V)$. The {\em rank} of $P$ is the dimension of the minimal by
inclusion face $F_P$ of $HYP(V)$ which contains $d_P$.

It is shown in \cite{DL} that $d_P$ determines uniquely the Delaunay
polytope $P$. When we move $d_P$ inside $F_P$, the Delaunay polytope $P$
changes, while its affine type remain the same.
In other words, like the rank of $P$, the affine type of $P$ is
an invariant of the face $F_P$.

The above movement of $d_P$ inside $F_P$ corresponds to a perturbation
of each basis of $L$, and, therefore, of each Gram matrix (i.e., each 
quadratic form) related to $L$. In this paper, we show that there is a 
one-to-one correspondence between the space spanned by $F_P$ and the space
${\cal B}(P)$ spanned by the set of perturbed quadratic forms. 
Hence, those two spaces have the same dimension.
It is shown here, that if one knows the coordinates of
vertices of $P$ in a basis, then it is simpler
to compute ${\rm dim}({\cal B}(P))$ than ${\rm dim}(F_P)$.
This fact is illustrated by
computations of ranks of cross polytope and half-cube.

%For all $u,v\in V(P)$, the vector $u-v$ connecting the vertices $u$ and
%$v$ is a minimal vector of the coset $L/2L$ which contains $u-v$.

%There are two ways to deform a Delaunay polytope $P$. On the one hand, we
%can move the corresponding hypermetric $d_P$ inside the face $F_P$ of the
%hypermetric cone $HYP$, where $d_P$ lies. On the other hand, we can
%perturbate the Gram matrix (i.e. the quadratic form $f_P$) of the lattice,
%which has $P$ as a Delanay polytope. In other words, we can move the
%quadratic form $f_P$ inside the L-type domain ${\cal D}_P$, where the
%form $f_P$ lies. In this paper, we show that the dimensions of the face
%$F_P$ and the L-type domain ${\cal D}_P$ are equal if $P$ is basic.

In the last section, we describe a non-basic repartitioning Delaunay
polytope recently discovered by the first author.

\section{Equalities of negative type and hypermetric}
A sphere $S=S(c,r)$ of radius $r$ and center $c$ in 
an $n$-dimensional lattice $L$ is said to be an {\em empty sphere}
if the following two conditions hold:
\begin{enumerate}
\item[(i)] $\Vert a-c\Vert^2\ge r^2$ for all $a\in L$,
\item[(ii)] the set $S\cap L$ contains $n+1$ affinely independent points.
\end{enumerate}
A Delaunay polytope $P$ in a lattice $L$ is a polytope, whose vertex-set
is $L\cap S(c,r)$ with $S(c,r)$ an empty sphere.

Denote by $L(P)$ the lattice generated by $P$. In this paper, 
we can suppose that $P$ is {\em generating} in $L$, i.e., that $L=L(P)$.
A subset $V\subseteq V(P)$ is said to be
{\em $\KK$-generating}, with $\KK$ being a ring,
if every vertex $w\in V(P)$ has a representation
$w=\sum_{v\in V}z(v)v$ with $1=\sum_{v\in V} z(v)$ and $z(v)\in \KK$.
If $|V|=n+1$, then $V$ is called an $\KK$-{\em affine basis};
the Delaunay polytope $P$ is called {\em $\KK$-basic} if it admits
at least one $\KK$-affine basis.
In this work $\KK$ will be $\ZZ$, $\QQ$
or $\RR$ and if the ring is not precised, it is $\ZZ$.
Furthermore, let
\begin{equation}\label{dep}
Y(P)=\{y\in \ZZ^{V(P)}: \sum_{v\in V(P)}y(v)v=0, \mbox{  }
\sum_{v\in V(P)}y(v)=0\}
\end{equation}
be the $\ZZ$-module of all integral dependencies on $V(P)$.
If the Delaunay polytope $P$ is a simplex, then $Y(P)=\{0\}$.

A dependency on $V(P)$ implies some dependencies between distances
$d_P(u,v)$ as follows. Let $c$ be the center of the empty sphere $S$
circumscribing $P$. Then all vectors $v-c$, $v\in V(P)$, have the
same norm $\Vert v-c\Vert^2=r^2$, where $r$ is the radius of the sphere $S$. Hence, 
\begin{equation}\label{dist}
d_P(u,v)=\Vert u-v\Vert^2=\Vert u-c-(v-c)\Vert^2=2(r^2-\langle u-c, v-c\rangle).
\end{equation}
Multiplying this equality by $y(v)$ and summing over $v\in V(P)$, we get
\begin{equation*}
\sum_{v\in V(P)}y(v)d_P(u,v)=2r^2\sum_{v\in V(P)}y(v)-
2\langle u-c, \sum_{v\in V(P)}y(v)(v-c)\rangle.
\end{equation*}
Since $y\in Y(P)$, we obtain the following important equality
\begin{equation}\label{yd}
\sum_{v\in V(P)}y(v)d_P(u,v)=0,\mbox{~for~any~}u\in V(P)\mbox{~and~}y\in Y(P).
\end{equation}
Denote by ${\cal S}_{dist}(P)$ the system of equations (\ref{yd}) for
all integral dependencies $y\in Y(P)$ and all $u\in V(P)$, where the
distances $d_P(u,v)$ are considered as unknowns.

Multiplying the equality (\ref{yd}) by $y(u)$ and summing over all
$u\in V(P)$, we obtain
\begin{equation}\label{neg}
\sum_{u,v\in V(P)}y(u)y(v)d_P(u,v)=0.
\end{equation}
This equality is called an equality of {\em negative type} and the system 
of such equality is denoted ${\cal S}_{neg}(P)$. Hence, the equalities of
${\cal S}_{neg}(P)$ are implied by the one of ${\cal S}_{dist}(P)$.

Each integral dependency $y\in Y(P)$ determines the following
representation of a vertex $w\in V(P)$ as an integer combination of
vertices from $V(P)$:
\[w=w+\sum_{v\in V(P)}y(v)v=\sum_{v\in V(P)}y^w(v)v, \]
where 
\begin{equation*}
y^w(v)=\left\lbrace\begin{array}{cl}
y(v)&\mbox{~if~}v\not=w,\\
y(w)+1&\mbox{~if~}v=w
\end{array}\right.
\mbox{~~and~~}\sum_{v\in V(P)}y^w(v)=1. 
\end{equation*}
Let $\delta_w$ be the indicator function of
$V(P)$: $\delta_w(v)=0$ if $v\not=w$, and $\delta_w(w)=1$. Obviously,
$\delta_w$ is $y^w$ for the trivial representation $w=w$. We have
$y^w=y+\delta_w$. Conversely, every representation
$w=\sum_{v\in V(P)}y^w(v)v$ provides the dependency
%$0=\sum_{v\in V(P)}y^w(v)v-w=\sum_{v\in V(P)}y(v)v$, where
$y=y^w-\delta_w\in Y(P)$. Substituting $y=y^w-\delta_w$ in (\ref{neg}),
we obtain the following equality
\[\sum_{u,v\in V(P)}y(u)y(v)d_P(u,v)=\sum_{u,v\in V(P)}y^w(u)y^w(v)d_P(u,v)-
2\sum_{v\in V(P)}y^w(v)d_P(w,v). \]
Since $d_P(w,w)=0$, we can set $y^w=y$ in the last sum. For any $w\in V(P)$,
we use this equality in the following form using equations (\ref{yd}) and
(\ref{neg})
\begin{equation}
\label{nh}
\sum_{u,v\in V(P)}y^w(u)y^w(v)d_P(u,v)=\sum_{u,v\in V(P)}y(u)y(v)d_P(u,v)+
2\sum_{v\in V(P)}y(v)d_P(w,v)=0.
\end{equation}
The equality
\[\sum_{u,v\in V(P)}z(u)z(v)d_P(u,v)=0, \mbox{ where }
\sum_{v\in V(P)}z(v)=1, \mbox{  }z(v)\in\ZZ, \]
is the hypermetric equality. Denote by ${\cal S}_{hyp}(P)$ the
system of all hypermetric equalities which hold for $d_P(u,v)$, considering
the distances $d_P(u,v)$ as unknowns.

In \cite{DL}, the following lemma is proved. For the sake of completeness,
we give its short proof.
\begin{lem}\label{hy}
Let $P$ be a Delaunay polytope with vertex-set $V(P)$. Let
$y^w\in \ZZ^{V(P)}$, such that $\sum_{v\in V(P)}y^w(v)=1$.
Then the following assertions are equivalent

\begin{enumerate}
\item[(i)] a vertex $w\in V(P)$ has the representation
$w=\sum_{v\in V(P)}y^w(v)v$;

\item[(ii)] the distance $d_P$ satisfies the hypermetric equality
$\sum_{u,v\in V(P)}y^w(u)y^w(v)d_P(u,v)=0$.
\end{enumerate}
\end{lem}
\proof (i)$\Rightarrow$(ii) Obviously, $y=y^w-\delta_w$, is a
dependency, i.e. $y\in Y(P)$. Hence, this implication follows from the
equalities (\ref{nh}), (\ref{yd}) and (\ref{neg}).

(ii)$\Rightarrow$(i) Substituting the expression (\ref{dist}) for $d_P$
in the hypermetric equality of (ii) we obtain the equality
\[2r^2-2\Vert \sum_{v\in V(P)}y^w(v)(v-c)\Vert^2=0. \]
Obviously, $\sum_{v\in V(P)}y^w(v)c=c$ and $\sum_{v\in V(P)}y^w(v)v$
is a point of $L(P)$. Denote this point by $w$. Then the above equality
takes the form $\Vert w-c\Vert^2=r^2$. Hence, $w$ lies on the empty sphere
circumscribing $P$. Therefore, $w\in V(P)$ and (i) follows. \qed

\vspace{2mm}
According to Lemma~\ref{hy}, each hypermetric equality of the system
${\cal S}_{hyp}(P)$ corresponds to a representation $y^w$ of a vertex
$w\in V(P)$. Since the relation $y=y^w-\delta_w$ gives a one-to-one
correspondence between dependencies on $V(P)$ and non-trivial
representations $y^w$ of vertices $w\in V(P)$, we can prove the
following assertion:
\begin{lem}
\label{equi}
The systems of equations ${\cal S}_{dist}(P)$ and ${\cal S}_{hyp}(P)$ are
equivalent, i.e., their solution sets coincide.
\end{lem}
\proof The equality (\ref{nh}) shows that each equation of the
system ${\cal S}_{hyp}(P)$ is implied by equations of the system
${\cal S}_{dist}(P)$.

Now, we show the converse implication. Suppose the unknowns $d(u,v)$
satisfy all hy\-per\-met\-ric equalities of the system ${\cal S}_{hyp}(P)$.
The equality (\ref{nh}) implies the equality
\[2\sum_{v\in V(P)}y(v)d(w,v)=-\sum_{u,v\in V(P)}y(u)y(v)d(u,v), \]
where $y=y^w-\delta_w$. This shows that, for the dependency $y$ on
$V(P)$, $\sum_{v\in V(P)}y(v)d(w,v)$ does not depend on $w$;
denote it by $A(y)$. Hence, we have
\[-2A(y)=\sum_{u,v\in V(P)}y(u)y(v)d(u,v)=A(y)\sum_{u\in V(P)}y(u). \]
According to equation (\ref{dep}), the last sum equals zero.
This implies the equalities (\ref{neg}) and hence 
the equalities of the system ${\cal S}_{dist}(P)$.  \qed

Obviously, the space determined by the system
${\cal S}_{hyp}(P)$ (and also of the system ${\cal S}_{dist}(P)$) is a
subspace $X(P)$ of the space spanned by all distances $d(u,v)$,
$u,v\in V(P)$. The dimension of $X(P)$ is the rank of $P$.
According to Lemma~\ref{equi}, in order to compute the rank of $P$,
we can use only equations of the system ${\cal S}_{dist}(P)$.

%The set $Y(P)$ of all integral dependencies between vertices of $P$
%forms a module over $\ZZ$. 
%This module is generated by special dependencies defined below.
Let $V_0=\{v_0,v_1,\dots,v_n\}$ be an $\RR$-affine basis of $P$.
Then each vertex $w\in V(P)$ has a unique representation
through vertices of $V_0$ as follows
\[w=\sum_{v\in V_0}x(v)v, \mbox{  }\sum_{v\in V_0}x(v)=1, \mbox{  }
x(v)\in \RR. \]
Since the vertices of $P$ are points of a lattice, in fact, $x(v)\in \QQ$.
Hence, the above equation can be rewritten as an integer dependency
\begin{equation*}
y_w(w)w+\sum_{v\in V_0}y_w(v)v=0, \mbox{  }
y_w(w)+\sum_{v\in V_0}y_w(v)=0, \mbox{~with~}y_w(v)\in\ZZ.
\end{equation*}
One sets $y_w(u)=0$ for $u\in V(P)-(V_0\cup\{w\})$ and gets $y_w\in Y(P)$.
Any dependency $y\in Y(P)$ is a rational combination of dependencies $y_w$,
$w\in V(P)-V_0$. Hence, the following equality holds:
\begin{equation*}
\beta y=\sum_{w\in V(P)-V_0}\beta_wy_w, \mbox{~with~}\beta_w\in \ZZ
\mbox{~and~}0<\beta\in \ZZ
\end{equation*}
Since the equalities (\ref{yd}) are linear over $y\in Y(P)$, the
dependencies $y_w$, $w\in V(P)-V_0$ provide the following system, which
is equivalent to ${\cal S}_{dist}(P)$
\begin{equation}\label{uwd}
y_w(w)d_P(u,w)+\sum_{v\in V_0}y_w(v)d_P(u,v)=0, \mbox{~with~}u\in V(P)
\mbox{~and~}w\in V(P)-V_0.
\end{equation}
We see that, for $u\in V(P)-V_0$, the distance $d_P(u,w)$, $w\in V(P)-V_0$,
is also expressed through distances between $u$ and $v\in V_0$.
But for $u\in V_0$, the distance $d_P(u,w)$ is expressed through
distances between $u,v\in V_0$.
This implies that the distance $d_P(u,w)$ for $u,w\in V(P)-V_0$
can be also represented through distances $d_P(u,v)$
for $u,v\in V_0$.
Hence, the dimension of $X(P)$ does not exceed $\frac{n(n+1)}{2}$,
where $n+1=|V_0|$, which is the dimension of the space of distances
between the vertices of $V_0$.

In order to obtain dependencies between $d_P(u,v)$ for $u,v\in V_0$, we
use equation (\ref{uwd}) for $u=w$. Since $d_P(w,w)=0$, we obtain
the equations
\[\sum_{v\in V_0}y_w(v)d_P(v,w)=0, \mbox{  }w\in V(P)-V_0. \]
%We change here the summation over $v\in V_0$ by the summation over $u\in V_0$.
Multiplying the above equation by $y_w(w)$ and using equation (\ref{uwd}),
we obtain
\[0=\sum_{u\in V_0}y_w(u)(y_w(w)d_P(u,w))=
-\sum_{u\in V_0}y_w(u)\sum_{v\in V_0}y_w(v)d_P(u,v). \]

So, we obtain the following main equations for dependencies between
$d_P(u,v)$ for $u,v\in V_0$
\begin{equation}
\label{duv}
\sum_{u,v\in V_0}y_w(u)y_w(v)d_P(u,v)=0, \mbox{  }w\in V(P)-V_0.
\end{equation}
Note, that if $V_0$ is an affine basis of $L(P)$, then one can set
$y_w(w)=-1$. In this case, the equation $y_w(w)+\sum_{v\in V_0}y_w(v)=0$
takes the form $\sum_{v\in V_0}y_w(v)=1$. This implies that the above
equations are hypermetric equalities for a $\ZZ$-basic Delaunay polytope $P$.
If $P$ is $\ZZ$-basic, then the distance $d_P$
restricted to the set $V_0$ lies on the face of the cone $HYP(V_0)$ 
determined by the hypermetric equalities (\ref{duv}). But if $P$ is not 
$\ZZ$-basic, then the equations (\ref{duv}) are not hypermetric,
and the distance $d_P$ restricted on the set $V_0$ lies inside
the cone $HYP(V_0)$. 
On the other hand, the distance $d_P$ on the whole set $V(P)$ lies
on the boundary of the cone $HYP(V(P))$.
This implies that, in this case, the rank of $d_P$ restricted to $V_0$
is greater than the rank of $d_P$ on $V(P)$.

This can be explained as follows. We can consider the cone $HYP(V_0)$
as a projection of $HYP(V(P))$ on a face of the positive orthant
$\RR^N_+$, where $N=|V(P)|$. This face is  determined by the
equations $d(u,v)=0$ for $v\in V(P)-V_0$ or/and $u\in V(P)-V_0$. By
this projection, the distance $d_P$, lying on the boundary of the cone
$HYP(V(P))$, is projected into the interior of the cone $HYP(V_0)$. This
hypermetric space corresponds to a wall of an $L$-type domain, which lies
inside the cone $HYP(V_0)$.

But, in order to compute the rank of $P$, it is sufficient to find the
dimension of the space determined by the system (\ref{duv}).

\section{Dependencies between lattice vectors}
Now we go from affine realizations to linear realizations. Take
$v_0\in V_0$ as origin of the lattice $L(P)$ and choose the lattice
vectors $a_i=a(v_i)=v_i-v_0$, $1\le i\le n$ such that 
$\{a_i:1\le i\le n\}$ forms a $\QQ$-basis of $L(P)$.
If $P$ is basic, we can choose $v_i$ such that $\{a_i:1\le i\le n\}$ is
a $\ZZ$-basis of $L(P)$.
Using the expressions $d_P(v_i,v_j)=\Vert a_i-a_j\Vert^2$, it is easy
to verify that there is the following relation between distances
$d_P(u,v)$, $u,v\in V_0$, and inner products $\langle a_i, a_j\rangle $:
\[d_P(v_i,v_0)=\Vert a_i\Vert^2, \mbox{    }d_P(v_i,v_j)=\Vert a_i\Vert^2
-2\langle a_i, a_j\rangle+\Vert a_j\Vert^2. \] 
And conversely, 
\[ \Vert a_i\Vert^2=d_P(v_i,v_0), \mbox{    }
\langle a_i, a_j\rangle=\frac{1}{2}(d_P(v_i,v_0)+d_P(v_j,v_0)-d_P(v_i,v_j)). \]
This shows that there is a one-to-one correspondence between the set
of distances $d_P(v_i,v_j)$, $0\le i<j\le n$, and the set of inner
products $\langle a_i, a_j\rangle$, $1\le i\le j\le n$.

We substitute the above expressions for $d_P(v_i,v_j)$, $0\le i,j\le n$,
into the equations (\ref{duv}), where we set $y_w(i)=y_w(v_i)$, and use
the equality $\sum_{i=0}^ny_w(i)=-y_w(w)$. We obtain the following
important equations
\begin{equation}\label{ya}
-y_w(w)\sum_{i=1}^ny_w(i)\Vert a_i\Vert^2=\Vert \sum_{i=1}^ny_w(i)a_i\Vert^2, \mbox{  }
w\in V(P)-V_0.
\end{equation}

We can obtain the equation (\ref{ya}) directly, as follows.
For $v\in V(P)$, the vector $a(v)=v-v_0$ is a lattice vector of 
$L(P)$. For $y\in Y(P)$, we have obviously $\sum_{v\in V(P)}y(v)a(v)=0$. In
particular, for $y=y_w$, this equation has the form
\[y_w(w)a(w)+\sum_{i=1}^ny_w(i)a_i=0 \]
and allows to represent the vectors $a(w)$ in the $\QQ$-basis
$\{a_i:1\le i\le n\}$.

Recall that the lattice vector $a(w)$ of each vertex $w\in V(P)$ of a
Delaunay polytope $P$ satisfies the equation $\Vert a(w)-c\Vert^2=r^2$.
Since $v_0\in V$, $a(v_0)=0$, which implies
$\Vert c\Vert^2=\Vert 0-c\Vert^2=r^2$.
The vertex-set of $P$ provides the following system of
equations $\Vert a(w)-c\Vert^2=\Vert c\Vert^2$, $w\in V(P)$, i.e.,
\begin{equation}\label{caw}
2\langle c, a(w)\rangle=\Vert a(w)\Vert^2, \mbox{  }w\in V(P).
\end{equation}
Since $y_w(w)a(w)=-\sum_{i=1}^ny_w(i)a_i$, the above equations take the form
\begin{equation}\label{ca}
-y_w(w)\sum_{i=1}^ny_w(i)2\langle c, a_i\rangle=\Vert \sum_{i=1}^ny_w(i)a_i\Vert^2.
\end{equation}
Recall that $a_i=a(v_i)$. Hence, the vertices $v_i$ give
$2\langle c, a_i\rangle=\Vert a_i\Vert^2$,
and the above equation takes the form of equation (\ref{ya}).

We will use the equations (\ref{ya}) mainly for basic Delaunay polytopes.
In this case, we can set $y_w(w)=-1$, and $a_i=b_i$, $1\le i\le n$,
where $B=\{b_i:1\le i\le n\}$ is the basis of $L(P)$ consisting of
lattice vectors of the basic Delaunay polytope $P$.

For a given $\QQ$-affine basis $V_0\subseteq V(P)$ of a Delaunay polytope $P$,
the set of affine dependencies $\{y_w\in Y(P):w\in V(P)-V_0\}$ is
uniquely determined up to integral multipliers and form a 
$\QQ$-basis of the $\ZZ$-module $Y(P)$.
This implies that the equations (\ref{ya}) determine a subspace
\begin{equation*}
{\cal A}(P)=\{a_{ij}: -y_w(w)\sum_{i=1}^ny_w(i)a_{ii}=
\sum_{1\le i,j\le n}y_w(i)y_w(j)a_{ij}, \mbox{  }y_w\in Y(P),
w\in V(P)-V_0\}
\end{equation*}
in the $\frac{n(n+1)}{2}$-dimensional space of all symmetric
$n\times n$-matrices $a_{ij}=a_{ji}$, $1\le i<j\le n$,

Since there is a one-to-one correspondence between distances $d(v_i,v_j)$,
$0\le i<j\le n$, and inner products $a_{ij}=\langle a_i, a_j\rangle$,
$1\le i\le j\le n$, the dimension of the subspace ${\cal A}(P)$ is equal to
the rank of $P$.
So, in order to compute the rank of $P$, we have to find
the dimension of ${\cal A}(P)$.

\section{The space ${\cal B}(P)$ and our computational method}

Fix a basis $B=\{b_i:1\le i\le n\}$ of the lattice $L$. Every 
lattice vector $a(v)$, $v\in V(P)$, has a unique representation
$a(v)=\sum_{i=1}^nz_i(v)b_i$. Define ${\cal Z}_B(P)=\{z_i(v):1\le i\le n, v\in V(P)\}$.

Recall that the cone ${\cal P}_n$ of all positive semi-definite forms on
$n$ variables is partitioned into $L$-type domains. Each $L$-type domain is
an open polyhedral cone of dimension $k$, where
$1\le k\le \frac{n(n+1)}{2}$. It consists of form having affinely
equivalent partitions into Delaunay polytopes, i.e. {\em Delaunay
partitions}. More exactly, an $L$-type domain is the set of quadratic forms
$f(x)=\Vert \sum_{i=1}^nx_ib_i\Vert^2$ having the same set of matrices
${\cal Z}_B(P)$ for all non-isomorphic Delaunay polytopes $P$ of its Delaunay
partition. So, this set is not changed when the basis $B$ changes such
that the form $f(x)$ belongs to the same $L$-type domain. In other words,
${\cal Z}_B(P)$ is an invariant of this $L$-type domain.

We set $z_{ij}=z_i(v_j)$ for $v_j\in V_0-\{v_0\}$, $1\le j\le n$. The
matrix $Z_B=(z_{ij})_1^n$ is non-degenerate and gives a correspondence
between the linear bases of $P$ and bases of $L(P)$. In particular, this
correspondence maps the space ${\cal A}(P)$ in the space ${\cal B}(P)$ of
matrices $b_{ij}=\langle b_i, b_j\rangle$ of the quadratic form $f(x)$.
If $P$ is basic and $b_i=a_i$, $1\le i\le n$, then $Z_B$
is the identity matrix $I$, and ${\cal A}(P)={\cal B}(P)$.

Substituting in the equations (\ref{caw}) the above representations of
the vectors $a(v)$, $v\in V(P)$, in the basis $B$, we obtain explicit
equations, determining the space ${\cal B}(P)$. In fact, we have
\begin{equation}\label{zb}
2\sum_{i=1}^nz_i(v)\langle c, b_i\rangle=\sum_{1\le i,j\le n}z_i(v)z_j(v)b_{ij},
\mbox{  }v\in V(P).
\end{equation}
We have the following $\frac{n(n+1)}{2}+n=\frac{n(n+3)}{2}$ parameters
in the equations (\ref{zb}):
\[b_{ij}=\langle b_i, b_j\rangle,\mbox{  }1\le i\le j\le n, \mbox{ and }\langle c, b_i\rangle, \mbox{  }
1\le i\le n, \]
Hence, all these parameters can be represented through a number of
independent parameters. This number is just the rank of $P$. Recall that a
Delaunay polytope is called {\em extreme} if ${\rm rk}(P)=1$. Hence,
in order to be extreme, a Delaunay polytope should have at least 
$\frac{n(n+3)}{2}$ vertices. 
 
Note that, for $v=v_0$, the equation (\ref{zb}) is an identity, since
$a(v_0)=0$ and therefore $z_i(v_0)=0$ for all $i$. So, we have $|V(P)|-1$
equations (\ref{zb}). For $v=v_i$, $1\le i\le n$, one gets $n$ equations
that give a representation of the parameters $\langle c, b_i\rangle$,
$1\le i\le n$ in terms of the parameters $\langle b_i, b_j\rangle$,
$1\le i\le j\le n$.
Hence, the equations (\ref{zb}), for $v\in V(P)-V_0$, allow
to find dependencies between the main parameters $\langle b_i, b_j\rangle$.
Now, we write out explicitly dependencies between $\langle b_i, b_j\rangle$.

Since the basic vectors $b_i\in B$ are mutually independent, a dependency
$\sum_{v\in V}y(v)a(v)=0$ implies the dependencies
$\sum_{v\in V}y(v)z_i(v)=0$ between the coordinates $z_i(v)$, $1\le i\le n$.

Multiplying equation (\ref{zb}) by $y(v)$, and summing over all
$v\in V(P)$, we obtain that the $\ZZ$-module $Y(P)$ determines the following
subspace of the space of parameters $b_{ij}=\langle b_i, b_j\rangle$:
\[ {\cal B}(P)=\{b_{ij}:\sum_{i,j=1}^n(\sum_{v\in V}y(v)
z_i(v)z_j(v))b_{ij}=0, \mbox{  }y\in Y(P)\}. \]
In the Delaunay partition of the lattice $L(P)$, there are infinitely
many Delaunay polytopes equivalent to $P$. Each of them has the form
$a\pm P$, where $a=\sum_{i=1}^nz_i^ab_i$ is an arbitrary lattice vector
of $L(P)$. Now, we show that the space ${\cal B}(P)$ is independent on a
representative of $P$ in $L(P)$, i.e., that ${\cal B}(P)={\cal B}(a\pm P)$.

Let $v_a=a\pm v$ be the vertex of the polytope $a\pm P$ corresponding to
a vertex $v$ of $P$. Obviously, $z_i(v_a)=z_i^a\pm z_i(v)$. Substituting
these values of $z_i(v_a)$ into the equations determining
${\cal B}(a\pm P)$, we obtain
\[\sum_{v_a}y(v_a)z_i(v_a)z_j(v_a)=\sum_{v\in V(P)}y(v)(z_i^az_j^a\pm
z_i^az_j(v)+z_i(v)z_j(v)). \]
Since $y$ is a dependency between vertices of $P$, the sums
with $z_1^a$ equal zero.
This shows that ${\cal B}(P)$ does not depend on a representative of $P$.
%This equation shows that the lattice $L$ having $P_0$ as a Delaunay
%polytope lies in a wall separating some $L$-type domains.

Since the equalities determining the space ${\cal B}(P)$ are linear in
$y$, we can consider these equalities only for basic dependencies $y_w$,
$w\in V(P)-V_0$. We obtain the following main system of equations
describing dependencies between the parameters $b_{ij}$;
\begin{equation}
\label{ywz}
 \sum_{i,j=1}^n(\sum_{v\in V}y_w(v)z_i(v)z_j(v))b_{ij}=0,
\mbox{  } w\in V(P)-V_0.
\end{equation}

A unimodular transformation maps a basis of $L(P)$ into another basis.
This trans\-for\-ma\-ti\-on generates a transformation which maps the space
${\cal B}(P)$ into another space related to $P$. The dimension of the
space ${\cal B}(P)$ is an invariant of the lattice $L(P)$
generated by $P$.

In \cite{BG}, a {\em non-rigidity degree} of a lattice was defined.
In terms of this paper, the non-rigidity degree of a lattice $L$ is
the dimension of the intersection of spaces ${\cal B}(P)$ related to
all non-isomorphic Delaunay polytopes of a star of Delaunay polytopes  of 
$L$. Hence,
\[{\rm nrd}(L)={\rm dim}(\cap_P{\cal B}(P)). \]
In fact, the space $\cap_P{\cal B}(P)$ is the supporting space of the
$L$-type domain of the lattice $L$.

\section{Centrally symmetric construction}
In many cases, the computation of the rank of a Delaunay polytope $P$
using the equations (\ref{zb}) is easier than by using the hypermetric
equalities generated by $P$.
We demonstrate this by giving a simpler proof of Lemma 15.3.7 of \cite{DL}.
Recall that a Delaunay polytope is either centrally
symmetric or asymmetric. Let $c$ be the center of the empty sphere
circumscribing $P$.
For any $v\in V(P)$, the point $v^*=2c-v$ is centrally symmetric to $v$.
If $P$ is centrally symmetric, then $v^*\in V(P)$ for all $v\in V(P)$.
If $P$ is asymmetric, then $v^*\not\in V(P)$ for all $v\in V(P)$.
\begin{lem}
\label{rkP}
Let $P$ be an $n$-dimensional basic centrally symmetric Delaunay polytope
of a lattice $L$ with the following properties:
\begin{enumerate}
\item The origin $0\in V(P)$ and the vectors $e_i$, $1\le i\le n$, are basic vectors of $L$, whose endpoints are vertices of $P$.
\item The intersection $P_1=P\cap H$ of $P$ with the hyperplane $H$
generated by the vectors $e_i$, $1\le i\le n-1$, is an
asymmetric Delaunay polytope of the lattice $L_1=L\cap H$.
\item If the endpoint $v_n$ of the basic vector $e_n$ is $v^*$ for
some $v\in V(P_1)$, then there is a vertex $u\in V(P)$ such that
$u\not=v,v^*$ for all $v\in V(P_1)$.
\end{enumerate}
Then ${\rm rk}(P) \le {\rm rk}(P_1)$.
\end{lem}
\proof It is sufficient to prove that the $n$ parameters
$\langle e_i, e_n\rangle$, $1\le i \le n$, can be expressed through the
parameters $a_{ij}=\langle e_i, e_j\rangle$, $1\le i\le j\le n-1$.

Let $c$ be the center of $P$. Obviously, $2c=0^*\in V(P)$. Since $P_1$
is asymmetric, $2c\not\in L_1$. It is easy to see that $2c=a_0+ze_n$,
with $a_0=\sum_{i=1}^{n-1}y_ie_i\in L_1$ and $0\not=z\in \ZZ$.
Hence, the equation $2\langle c, e_i\rangle=\Vert e_i\Vert^2$
takes the form $\langle a_0+ze_n, e_i\rangle=\Vert e_i\Vert^2$, and
the parameters $\langle e_i, e_n\rangle$ are represented through
the parameters $\langle e_i, e_j\rangle$ as follows
\[\langle e_i, e_n\rangle=\frac{1}{z}(\Vert e_i\Vert^2-\langle a_0, e_i\rangle), \mbox{  }1\le i\le n-1. \]
Now, using the equation $2\langle c, e_n\rangle=\Vert e_n\Vert^2$,
we obtain $\langle a_0+ze_n, e_n\rangle=\Vert e_n\Vert^2$,
i.e., $\Vert e_n\Vert^2(1-z)=\langle a_0, e_n\rangle=\sum_{i=1}^{n-1}y_i \langle e_i, e_n\rangle$.
Hence, if $z\not=1$, we can represent $\Vert e_n\Vert^2$ through
$\langle e_i, e_j\rangle$,
$1\le i\le j\le n-1$, too.
But if $z=1$, then the endpoint $v_n$ of $e_n$
belongs to $(V(P_1))^*$. In this case there is a vertex $u$ such that
$u=\sum_{i=1}^nz_ie_i=u_0+z_ne_n$, where $u_0\in L_1$ and $z_n\not=0,1$.
Using the equation $2\langle c, u\rangle=\Vert u\Vert^2$, where now $2c=a_0+e_n$, we have
$\langle a_0+e_n, u_0+z_ne_n\rangle=\Vert u_0+z_ne_n\Vert^2$.
This equation gives
\[\Vert e_n\Vert^2=\frac{1}{z_n(z_n-1)}[\langle a_0-u_0, u_0\rangle+\langle z_na_0+(1-2z_n)u_0, e_n\rangle].\] 
The strict inequality ${\rm rk}(P)<{\rm rk}(P_1)$ is possible if some vertices of 
the set $V(P)-V(P_1)$ provide additional relations between the parameters
$\langle e_i, e_j\rangle$, $1\le i\le j\le n-1$. \qed

\vspace{2mm}
%It is useful to compare this proof with the proof given in \cite{DL}. 
Examples, where ${\rm rk}(P)<{\rm rk}(P_1)$, can be given by some extreme
Delaunay polytopes. 
\begin{cor}
\label{ext}
Let $P$ be a basic centrally symmetric Delaunay polytope satisfying
the conditions of Lemma~\ref{rkP}. $P$ is extreme if $P_1$ is extreme.
\end{cor}

\section{Computing the rank of simplexes, cross-polytopes and half-cubes}

\indent {\bf Simplices}. Let $\Sigma$ be an $n$-dimensional simplex with vertices
0, $v_i$, $1\le i \le n$. The vertex $v_i$ is the end-point of the basic
vector $e_i$, $1\le i\le n$. We have only $n$ equations
$2\langle c, e_i\rangle=\Vert e_i\Vert^2$
determining only the coordinates of the center $c$ of $\Sigma$ in the
basis $\{e_i:1\le i\le n\}$. Since there is no relation between the
$\frac{n(n+1)}{2}$ parameters $\langle e_i, e_j\rangle=a_{ij}$,
all these parameters are independent.
Hence,
\[{\rm dim}({\cal B}(\Sigma))=\frac{n(n+1)}{2}, \mbox{ i.e., }
{\rm rk}(\Sigma)=\frac{n(n+1)}{2}. \]

{\bf Cross-polytopes}.
An $n$-dimensional cross-polytope $\beta_n$ is a basic centrally symmetric
Delaunay polytope. It is the convex hull of $2n$ end-points of $n$
linearly independent segments intersecting in the center $c$ of the
circumscribing sphere. The set $V(\beta_n)$ is partitioned
into two mutually centrally symmetric $n$-subsets each of which is
the vertex-set of an $(n-1)$-dimensional simplex $\Sigma$. So,
$V(\beta_n)=V(\Sigma)\cup V(\Sigma^*)$. Let
$V(\Sigma)=\{0,v_i:1\le i\le n-1\}$.
All $\ZZ$-affine bases of $\beta_n$ are of the same type: $n-1$ basic
vectors $e_i$, $1\le i\le n-1$, with end vertices $v_i$, are basic
vectors of the simplex $\Sigma$, and $e_n=2c$, which is the segment
which connects the vertex $0$ with its opposite vertex $0^*$.
Let $a_i$ be the lattice vector
endpoint of which is the vertex $v_i^*\in \Sigma^*$. Obviously,
$a_i=2c-e_i=e_n-e_i$. The equality $2\langle c, a_i\rangle=\Vert a_i\Vert^2$
gives $\langle e_i, e_n\rangle=\Vert e_i\Vert^2$,
$1\le i\le n-1$.
So, we obtain $n-1$ independent relations between the parameters
$\langle e_i, e_j\rangle$, and they are the only relations. Hence,
\[{\rm rk}(\beta_n)=\frac{n(n+1)}{2}-(n-1). \]
(Cf., the first formula on p.232 of \cite{DL}.)

\vspace{2mm}
{\bf Half-cubes}.
Take $N=\{1,2,\dots,n\}$, a basis $(e_i)_{i\in N}$ and defines
$e(T)=\sum_{i\in T}e_i$ for any $T\subseteq N$.
Call a set $T\subseteq N$
{\em even} if its cardinality $|T|$ is even. A half-cube $h\gamma_n$
is the convex hull of endpoints of all vectors $e(T)$ for all even
$T\subseteq N$. Note that $h\gamma_3$ is a simplex, and $h\gamma_4$ is the
cross-polytope $\beta_4$. Hence,
\[{\rm rk}(h\gamma_3)=\frac{3(3+1)}{2}=6, \mbox{~and~}
{\rm rk}(h\gamma_4)=\frac{4(4+1)}{2}-3=7. \]

The rank of $h\gamma_n$ is computed from the following system of equations:
\begin{equation}\label{eT}
2\langle c, e(T)\rangle=\Vert e(T)\Vert^2, \mbox{  }T\subseteq N, \mbox{  }T\mbox{ is even}.
\end{equation}
Let $T_1$ and $T_2$ be two disjoint even subsets of $N$. Since the set
$T=T_1\cup T_2$ is even, we have
\[2\langle c, e(T_1\cup T_2)\rangle =2\langle c, e(T_1)+e(T_2)\rangle
=\Vert e(T_1)+e(T_2)\Vert^2=
\Vert e(T_1)\Vert^2+\Vert e(T_2)\Vert^2+2\langle e(T_1), e(T_2)\rangle. \]
Comparing this equation with the equations (\ref{eT}) for $T=T_1$ and
$T=T_2$, we obtain that for any two disjoint even subsets the following
{\em orthogonality conditions} hold:
\begin{equation*}
\langle e(T_1), e(T_2)\rangle=0, \mbox{~if~}T_1\cap T_2=\emptyset,\mbox{  }
T_i\subset N, \mbox{~and~}T_i\mbox{~is~even,~ }i=1,2.
\end{equation*}
Note that, for $n=3$, we have no orthogonality condition. If $n\geq 4$,
take $4$ elements $i$, $j$, $k$ and $l$ and write three equalities
corresponding to three partitions:
\[\langle e_i+e_j, e_k+e_l\rangle=0, \mbox{  }\langle e_i+e_k, e_j+e_l\rangle=0, \mbox{  }
\langle e_i+e_l, e_j+e_k\rangle=0. \]
It is easy to verify that these equalities are equivalent to the
following three equalities
\begin{equation}\label{il}
\langle e_i, e_j\rangle+\langle e_k, e_l\rangle=0, \mbox{  }
\langle e_i, e_k\rangle+\langle e_j, e_l\rangle=0, \mbox{  }
\langle e_i, e_l\rangle+\langle e_j, e_k\rangle=0.
\end{equation}
In the particular case $n=4$, we conclude again that
${\rm rk}(h\gamma_4)=\frac{4(4+1)}{2}-3=7$.

We show that, for $n\ge 5$, the orthogonality conditions are equivalent
to mutual orthogonality of all vectors $e_i$, $i\in N$. To this end, it
is sufficient to consider even subsets of cardinality two and use
equation (\ref{il}) for each quadruple $\{i,j,k,l\}\subseteq N$.
Considering arbitrary subsets of $N$ of cardinality 4, we obtain that,
for $n\ge 5$, the system of equalities (\ref{il}) for all quadruples
has the following unique solution
\[\langle e_i, e_j\rangle=0, \mbox{  }1\le i<j\le n, \mbox{ for }n\ge 5. \]
So, all the basic vectors are mutually orthogonal. Obviously, the
orthogonality of basic vectors implies the orthogonality conditions.
Hence, the only independent parameters are the $n$ parameters
$\Vert e_i\Vert^2$, $i\in N$. This implies that
\[{\rm rk}(h\gamma_n)=n \mbox{~if~} n\ge 5. \] 
Note that we use a basis of $h\gamma_n$, which is not
a basis of the lattice generated by $h\gamma_n$.
%Note that we use here a basis which is not only the basis of $h\gamma_n$ 
%but it is not a basis of the lattice generated be $h\gamma_n$.
But the spaces ${\cal B}(P)$ have the same dimension for all bases.
See another proof in \cite{DutClass}.

%\newpage
\section{A non-basic repartitioning Delaunay polytope}
%The first author found an example of a non-basic Delaunay
%polytope $P_0$.
The example $P_0$ given in this section is $12$ dimensional;
its $14$ vertices belong to two disjoint sets of vertices
of regular simplexes $\Sigma_i^2$, $i=1,2$, of dimension $2$,
and two disjoint sets of vertices
of regular simplexes $\Sigma_i^3$, $i=1,2$, of dimension $3$.

Let $V(\Sigma_i^q)$ be the vertex-set of the four simplex $\Sigma_i^q$,
$i=1,2$, $q=2,3$. Then $V=\cup V(\Sigma_i^q)$ is the vertex-set of $P_0$.
The distances between the vertices of $P_0$ are as follows
\[d(u,v)= \left\{ \begin{array}{rllll}
7 & {\rm if} & u,v\in \Sigma_i^q, & i=1,2, & q=2,3; \\
6 & {\rm if} & u\in\Sigma_i^2, & v\in\Sigma_i^3,& i=1,2; \\
10 & {\rm if} & u\in\Sigma_1^2, & v\in\Sigma_2^2;&  \\
12 & {\rm if} & u\in\Sigma_1^3, v\in\Sigma_2^3 \mbox{ or }
& u\in\Sigma_1^2, v\in\Sigma_2^3, \mbox{ or}
& u\in\Sigma_2^2, v\in\Sigma_1^3.
\end{array}\right. \]
We show that, for every $u\in V$, the set $V-\{u\}$
is an $\RR$-affine basis of $P_0$. In fact, let $V-\{u\}=\{v_i:0\le i\le 12\}$
and let $a_i=v_i-v_0$, $1\le i\le 12$. For the Gram matrix
$a_{ij}=\langle a_i, a_j\rangle$,
we have $a_{ii}=\Vert a_i\Vert^2=\Vert v_i-v_0\Vert^2=d(v_i,v_0)$.
The relations between $a_{ij}$ and $d(v_i,v_j)$ are
$a_{ij}=\frac{1}{2}(d(v_i,v_0)+d(v_j,v_0)-d(v_i,v_j))$. Now, one can
verify that the Gram matrix $(a_{ij})$ is not singular. Hence,
$\{a_i:1\le i\le n\}$ is a basis, i.e. the dimension of $P_0$ is, in fact,
12.

The space $Y(P_0)$ of affine dependencies on vertices of $P_0$ is
one-dimensional.
%We have the following dependency generating $Y(P_0)$.
For $v\in V$, let
\[y(v)=\left\{ \begin{array}{rll}
3  & {\rm if} & v\in \Sigma_1^2,\\
-3 & {\rm if} & v\in \Sigma_2^2,\\
2  & {\rm if} & v\in \Sigma_2^3,\\
-2 & {\rm if} & v\in \Sigma_1^3.
\end{array}\right. \]
Obviously, $\sum_{v\in V}y(v)=0$.
It is easy to verify that for any $u\in V$ the following equality holds
\begin{equation}\label{bu}
\sum_{v\in V}y(v)d(u,v)=0.
\end{equation}
Let $S(c,r)$ be the sphere circumscribing $P_0$.
Then $\Vert v-c\Vert^2=r^2$ for all $v\in V$.
We have $d(u,v)=\Vert u-v\Vert^2=\Vert (u-c)-(v-c)\Vert^2=2(r^2-\langle u-c, v-c\rangle)$.
Since $\sum_{v\in V}y(v)=0$, the equality (\ref{bu}) takes the form
\begin{equation*}
\langle u-c, \sum_{v\in V}y(v)(v-c)\rangle=0, \mbox{ i.e., }
\langle u-c, \sum_{v\in V}y(v)(v-c)\rangle=\langle u-c, \sum_{v\in V}y(v)v\rangle=0.
\end{equation*}
Since this equality holds for all $u\in V$, and $\{u-c\mbox{~}|\mbox{~}u\in V\}$ 
span $\RR^{12}$, we obtain $\sum_{v\in V}y(v)v=0$, i.e., $y\in Y(P_0)$.
Since $Y(P_0)$ is one dimensional and the coefficient of $y$ have greatest
common divisor $1$, one has $Y(P_0)=y \ZZ$.

%have the following unique up to a multiplier affine
%dependency between vertices of $P_0$
%\begin{equation}\label{af}
%\sum_{v\in V}y(v)v=0.
%\end{equation}

Using the basis $\{a_i:1\le i\le 12\}$, for non-basic $a(w)$, we obtain
$a(w)=-\frac{1}{y(w)}\sum_{i=1}^{12}y(v_i)a_i$. Since there exist a $i$
such that $\frac{y(v_i)}{y(w)}\notin \ZZ$ for any choice of $w\in V$,
the polytope $P_0$ is not basic and the $\QQ$-basis $\{a_i:1\le i\le n\}$
is not a $\ZZ$-basis of any lattice $L$ having $P_0$ as a Delaunay polytope.

Remarking that we can put the vector $y$ in equation (\ref{ya}), we obtain
the following equation
\begin{equation*}
-y(w)\sum_{i=1}^{12}y(v_i)\Vert a_i\Vert^2=\Vert \sum_{i=1}^{12}y(v_i)a_i\Vert^2.
\end{equation*}
which implies that ${\rm rk}(P_0)={\rm rk}(V,d)=\frac{12\times 13}{2}-1=77$.

It is useful to compare the above computation of ${\rm rk}(P)$
with the following computations using distances.
Recall that ${\rm rk} (V(P_0),d)$ is equal to the
dimension of the face of the hypermetric cone $HYP(V(P_0))=HYP(V)$, where
the distance $d$ lies. The dimension of $HYP(V)$ is
$N=\frac{|V|(|V|-1)}{2}=\frac{14\cdot 13}{2}=91$.

As in Section 2, we obtain that, for every $w\in V=V(P_0)$, the equality
(\ref{bu}) implies the following hypermetric equality
\begin{equation}
\label{hp0}
\sum_{v,v'\in V}y^w(v)y^w(v')d(v,v')=0,
\end{equation}
where $y^w(v)=y(v)+\delta_w$.
%Each of the $14$ vertices $w$ of $P_0$ provides the hypermetric equality of the type (\ref{hp0}).
It is easy to see that the $14$ equalities (\ref{hp0})
for $14$ vertices $w\in V$ are mutually independent.
In fact, these $14$ equalities are equivalent to the
$14$ equalities (\ref{bu}) for the $14$ vertices $u\in V$. The two equations
(\ref{bu}) corresponding to two vertices $u,w\in V$ have only one common
distance $d(u,w)$. The intersection of the corresponding $14$ facets is a
face of dimension $91-14=77$.

But, for every $u\in V$, the hypermetric space $(V-\{u\},d)$ has rank
$\frac{(|V-\{u\}|)(|V-\{u\})|-1)}{2}=78$, which is greater
than ${\rm rk}(V,d)=77$.

%\vspace{3mm}
%\newpage

%\vspace{3mm}
%E-mail address: grishuhn@cemi.rssi.ru

\end{document}